\documentclass[a4paper,11.5pt,reqno]{amsart}

\usepackage{graphicx}
\usepackage[utf8]{inputenc}
\usepackage{amsmath,amssymb,amsthm}
\usepackage{enumitem}
\usepackage{xcolor}
\usepackage{tikz}             
\usepackage{hyperref}
\usepackage{caption}

\usetikzlibrary{shapes.geometric}

\usepackage{lineno}
\modulolinenumbers[3]

\theoremstyle{plain}

\theoremstyle{definition}
\newtheorem{definition}{Definition}

\theoremstyle{remark}
\newtheorem{remark}{Remark}

\newcommand{\R}{\mathbb{R}}
\newcommand{\Z}{\mathbb{Z}}

\numberwithin{equation}{section}

\hypersetup{colorlinks = true, linkcolor = blue, urlcolor = blue, citecolor = blue}

\renewcommand{\phi}{\varphi}

\begin{document}

\title[No new lower bound for the density of planar, 1-avoiding Sets]  {No new lower bound for the density of planar Sets avoiding Unit Distances}

\author[H. Ruhland]{Helmut Ruhland}
\address{Santa F\'{e}, La Habana, Cuba}
\email{helmut.ruhland50@web.de}

\subjclass[2020]{Primary 52C17, Secondary  52C10}

\keywords{distance avoiding sets, constant diameter, Croft's tortoise}

\begin{abstract}
In a recently published article by G. Ambrus et al. a new \emph{upper bound} for the density of an unit avoiding,  periodic set is given as $0.2470$, the first upper bound $< 1/4$. A  construction of Croft 1967 gave a \emph{lower bound}  $\delta_C = 0.22936$ for the density. To this date, no better construction with a higher bound has been given. In the \emph{first versions} of this article I gave a construction of planar sets with a "higher" density than Croft's tortoises. No explicit value for this density was given, it was  just shown that Croft's density is a local minima in the density of the constructed 1-parameter family of planar sets. But now I found a servere error. After the correction in this article none of the investigated sets of constant diameter resulted in a new lower bound. I did not withdraw the article, maybe something could be useful for somebody.
\end{abstract}

\date{\today}

\maketitle

\section{Introduction}

In a recently published article \cite{CmCsMaVaZs} a new \emph{upper bound} for the density $m_1 (\R^2)$ of an unit avoiding,  periodic set is given as $0.2470$, the first upper bound $< 1/4$. The existence of such a bound $< 1/4$ was already conjectured by Erd\"os. A  construction of Croft 1967 gave a \emph{lower bound}  $\delta_C = 0.22936$ for the density, i.e. $m_1 (\R^2) \ge \delta_c$. To this date, no better construction has been given. In the \emph{first versions} of this article, denote emphasis at "first", I gave a construction of planar sets with a "higher" density than Croft's tortoises. No explicit value for this density was given, it was  just shown that Croft's density is a local minima in the density of a here constructed 1-parameter family of planar sets.

But now I found a servere error in the former \ref{A2hat_d_4r} for 2-parameter minimization,  when I compared theoretical values with the values obtained by a numerical implementation. This implementation was intended to give a concrete value for the predicted density higher than Croft's. After correcting the error none of the here investigated sets of constant diameter results in a new lower bound. I did not withdraw the article, I just changed the title, maybe something could be useful for somebody. \\

A short outline of this article: \\
In the following section \ref{set_const_diameter} a 1-parameter family $D_\epsilon$ of  sets with constant diameter $2$ is defined. In section \ref{construction} with these sets $D_\epsilon$ a family $S_\epsilon$ of 2-avoiding sets is constructed. In section \ref{calc_area_tort} the area of the tortoises (the sets of constant diameter with 6 cutted segments) is calculated as power series up to quadratic terms. Unfortunatly this shows that for small $\epsilon \ne 0$ the density of $S_\epsilon$ is $<$ the density of Croft's tortoises  $S_0$. \\
In the appendices I present formulae for areas of disc segments as power series. In appendix \ref{3col_hex_lattice} some information about 3-colored hexagonal lattices is given.

\section{A 1-parameter family of sets with constant diameter $2$ \label{set_const_diameter}}

For the construction of sets with \emph{constant width} $2 \, r$, see \cite{KaSw2019}, Recipe 3, equation (16). A periodic function $q (\phi)$ with period $2 \, \pi$ defines the boundary of such a set. $q (\phi)$ has the following property: 
\begin{equation}
   q (\phi + \pi) = - q (\phi)
   \label{q_phi_prop}
\end{equation}

\begin{definition}
Sets or curves with \emph{constant diameter} are defined in the following manner. The diameter of a point $B$ on the boundary of the set is defined as supremum of all distances between $B$ and points in the set. If the diameters of all points on the boundary are equals, we speak from a set with constant diameter. Sets with constant diameter are extremal in the following sense: Adding a small area on the boundary increases the diameter of the set. The area takes a local maximum, the global maximum of the area is of course the area of a circle. So with this property to be extremal, such sets are good candidates for tortoises before cutting.
\end{definition}

If we take $q (x)$ piecewise constant, we get curves and sets with constant diameter. The radii at the interval boundaries are not continuous. Here 24 intervals of different lengths are chosen. The interval ends are not equidistant because the sizes of the intervals $1, 4, 5, 8, \dots$ have to be greater than the half angle $\phi_C = 15.08686^{\circ}$ of Croft's tortoise.  \\
The interval size of the intervals above and below multiples of $\pi / 3$ is $4 \pi / 45 = 16^{\circ}$. 
So the cut for Croft's tortoise with an half angle of $\phi_C$ lies completely in the two  intervals.
This is neccesary to simplify the calculation of the area of the cutted segment.

\begin{equation}
 q (\phi) = \left\{
\begin{array}{rcrl}
                     0  &   \le \phi <   &  4 \, \pi /     45 & - 0.9779019 \, \epsilon \quad 16^{\circ} \\
     4 \, \pi /     45  &   \le \phi <   &    \, \pi / \;\; 6 & + 0.7242098 \, \epsilon \\
       \, \pi / \;\; 6  &   \le \phi <   & 11 \, \pi /     45 & - 0.7335699 \, \epsilon \\
    11 \, \pi /     45  &   \le \phi <   &    \, \pi / \;\; 3 & + 1.0000000 \, \epsilon \quad 16^{\circ} \\
       \, \pi / \;\; 3  &   \le \phi <   & 19 \, \pi /     45 & - 0.9227438 \, \epsilon \quad 16^{\circ} \\
    19 \, \pi /     45  &   \le \phi <   &    \, \pi / \;\; 2 & + 0.4889208 \, \epsilon \\
       \, \pi / \;\; 2  &   \le \phi <   & 26 \, \pi /     45 & - 0.1260494 \, \epsilon \\
    26 \, \pi /     45  &   \le \phi <   &  2 \, \pi / \;\; 3 & + 0.0442648 \, \epsilon \quad 16^{\circ} \\
     2 \, \pi / \;\; 3  &   \le \phi <   & 34 \, \pi /     45 & + 0.0042024 \, \epsilon \quad 16^{\circ} \\
    34 \, \pi /     45  &   \le \phi <   &  5 \, \pi / \;\; 6 & - 0.1019359 \, \epsilon \\
     5 \, \pi / \;\; 6  &   \le \phi <   & 41 \, \pi /     45 & + 0.4722281 \, \epsilon \\
    41 \, \pi /     45  &   \le \phi <   &                \pi & - 0.8992978 \, \epsilon \quad 16^{\circ} \\
    \dots \\
    11 \, \pi / \;\; 6  &   \le \phi <   & 86 \, \pi /     45 & - 0.4722281 \, \epsilon \\
    86 \, \pi /     45  &   \le \phi <   &           2 \, \pi & + 0.8992978 \, \epsilon \quad 16^{\circ} \\
\end{array} 
\right.
   \label{q_phi}
\end{equation}

\begin{center}
  \includegraphics[width=0.82\textwidth]{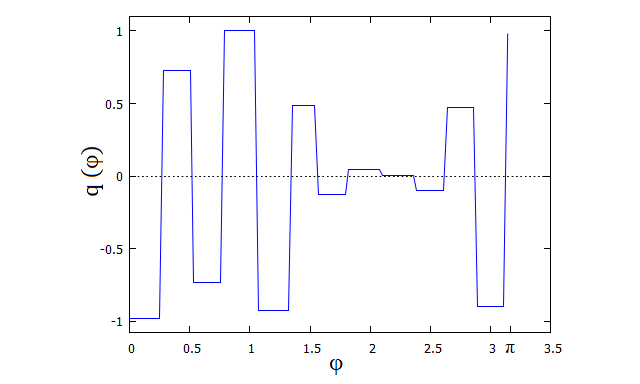}
  \captionof{figure}{$q (\phi)$ with $12$ intervals in $0 \dots \pi$ for $\epsilon = 1.0$.}
  \label{croft}
\end{center}

With formula (2) in \cite{KaSw2019} we get for the horizontal and vertical components of the vector
 $\mathbf{x} (\phi)$:
\begin{equation}
 x (\phi) = (1 - q (\phi)) \, \cos (\phi)  + \;  \left\{
\begin{array}{rcrl}
                     0  &   \le \phi <   &  4 \, \pi /     45 & - 0.9779019 \, \epsilon \quad 16^{\circ} \\
     4 \, \pi /     45  &   \le \phi <   &    \, \pi / \;\; 6 & + 0.6582729 \, \epsilon \\
    \dots \\
    11 \, \pi / \;\; 6  &   \le \phi <   & 86 \, \pi /     45 & - 0.4190975 \, \epsilon \\
    86 \, \pi /     45  &   \le \phi <   &           2 \, \pi & + 0.8992978 \, \epsilon \quad 16^{\circ} \\
\end{array} 
\right.
   \label{x_phi}
\end{equation}

\begin{equation}
 y (\phi) = (1 - q (\phi)) \, \sin (\phi)  + \;  \left\{
\begin{array}{rcrl}
                     0  &   \le \phi <   &  4 \, \pi /     45 & 0 \qquad  \qquad \qquad       16^{\circ} \\
     4 \, \pi /     45  &   \le \phi <   &    \, \pi / \;\; 6 & + 0.4691656 \, \epsilon \\
    \dots \\
    41 \, \pi /     45  &   \le \phi <   &             \, \pi & 0 \qquad \qquad \qquad        16^{\circ} \\
                 \, \pi &   \le \phi <   & 49 \, \pi /     45 & 0 \qquad \qquad \qquad        16^{\circ} \\
    \dots \\
    11 \, \pi / \;\; 6  &   \le \phi <   & 86 \, \pi /     45 & + 0.3780437 \, \epsilon \\
    86 \, \pi /     45  &   \le \phi <   &           2 \, \pi & 0 \qquad \qquad \qquad        16^{\circ} \\
\end{array} 
\right.
   \label{y_phi}
\end{equation}
These formulae show that in an interval the curve has a constant radius of $(1 - q (\phi))$.

\begin{center}
  \includegraphics[width=0.92\textwidth]{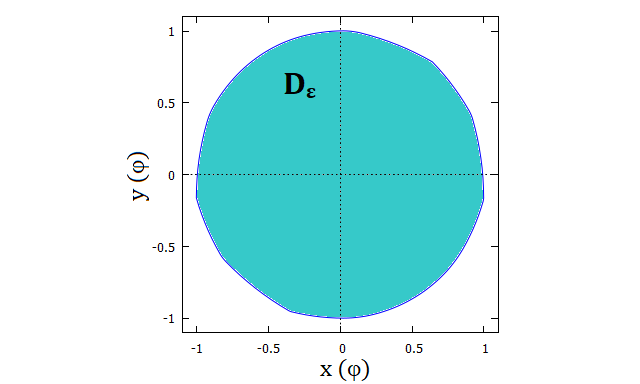}
  \captionof{figure}{The set $D_\epsilon$ with constant diameter $2$ for $\epsilon = 1.0$.}
  \label{croft}
\end{center}

The area of the set $D_\epsilon$ is:
\begin{equation}
   A_D = \frac{1}{2} \int\limits_{0}^{2 \pi} \left( x (\phi) \, \frac{dy (\phi)}{d\phi} 
                     - y (\phi) \, \frac{dx (\phi)}{d\phi}  \right) d\phi
              = \pi - 0.0104747054 \, \epsilon^2
   \label{xy_area}
\end{equation}
The area has no linear term in $\epsilon$, because the circle with $\epsilon = 0$ and area $\pi$
is extremal, i.e. has the largest area among all constant diameter $2$ sets. \\

Determining $q (\phi)$: \\
The terms with $\epsilon^2$ in the area of $D_\epsilon$, see \ref{xy_area} and in the area of a tortoise,
see \ref{area_tortoise}, represent quadratic forms in 1 variable, here in $\epsilon$. Here a general  "ansatz" for $q (\phi)$ with $12$ variables for the $24$ intervals of $q (\phi)$ because of property \ref{q_phi_prop}, $2$ variables for the shifts in x/y-direction. Because of equation (14) in \cite{KaSw2019}, $2$ variables are not free. So in the area of a tortoise, see \ref{area_tortoise}, we get a quadratic form in $12$ variables. We want this quadratic form to be positive for certain values. The quadratic form can be represented by a symmetric $12 \times 12$ matrix. This matrix has $1$ positive eigenvalue and $11$ negative eigenvalues. The positive eigenvalue and the corresponding eigenvector lead to the constants in \ref{q_phi}.\\

A visualization of the construction for curves of constant diameter:

\begin{center}
  \includegraphics[width=0.60\textwidth]{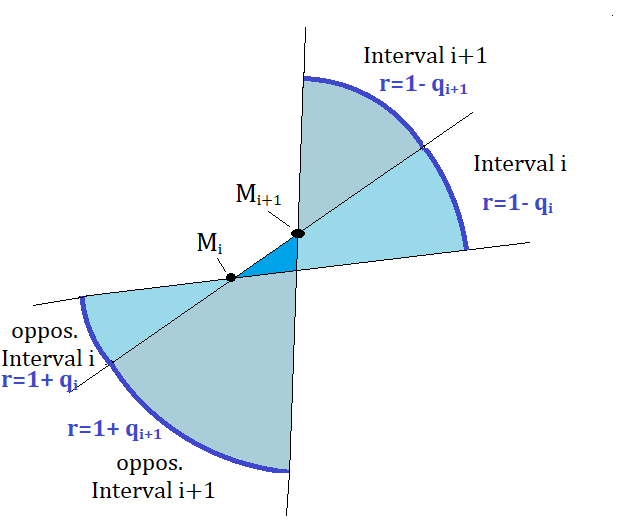}
  \captionof{figure}{$M_i$ is the midpoint of the $2$ circle arcs in interval $i$ and the opposite interval.
                     $M_{i+1}$ is the midpoint of the $2$ circle arcs in interval $i+1$ and the opposite
                     interval. Its easy to see: the diameter in each interval is $2$.}
  \label{croft}
\end{center}
The x/y coordinates of the midpoints $M_{i} = (x_{m,i}, y_{m,i})$ and $M_{i+1} = (x_{m,i+1}, y_{m,i+1})$ are given with the right, piecewise, not depending on $\phi$ part  of \ref{x_phi} and \ref{y_phi}. To be countinuous at the common interval ends $\phi_{i,i+1}$ in the middle of the $2$ intervals, the following conditions have to be fulfilled:
\begin{equation}
   \begin{split}
   & (x_{m,i+1} - x_{m,i}) - (q_{i+1} - q_{i}) \, \cos (\phi_{i,i+1}) = 0  \\
   & (y_{m,i+1} - y_{m,i}) - (q_{i+1} - q_{i}) \, \cos (\phi_{i,i+1}) = 0
   \label{cont_cond}
   \end{split}
\end{equation}

The complete data neccesary to represent the 1-parameter family $D_\epsilon$ of sets in an arbitrary
Computer Algebra System CAS is given in appendix \ref{copa_1_par}. By just copy \& paste, the numerical values of the description can be inserted in a CAS.

\section{Construction of a 1-parameter family of 2-avoiding sets \label{construction}}

In the previous section a 1-parameter family $D_\epsilon$ of sets with constant diameter $2$ was constructed. For $\epsilon = 0$ this set is a disc. A copy of these open sets is placed on the red lattice points of a colored hexagonal lattice, see figure \ref{3-Coloring}. A by $2 \, \pi / 3$ rotated copy at the green points (important green points!) and  a by $4 \, \pi /3$ rotated copy at the blue points. See also figure \ref{lattice_orientations} with the orientation of the copies.
The lattice constant of the hexagonal lattice is $L = 2 \, L_C$, $L_C$ the lattice constant used in Croft's unit distance avoiding set. So until now these so constructed periodic sets are not 2-avoiding, because
$L = 2 \, (1 + \cos (\phi_C)) < 4$.

Now a stripe of width $2$ is positioned between two rotated $D_\epsilon$ at distance $L$ (nearest neighbors), see figure \ref{gseg2}. The position of this stripe is varied until the sum of the areas of the $2$ disc segments, cutted by the stripe is minimal. This family of 2-avoiding sets is named $S_\epsilon$. 

For $\epsilon = 0$ $D_\epsilon$ is a disc, rotation, does not play any role. $S_0$ is the set obtained by Croft's construction stretched by a factor of $2$.

\section{Calculating the areas of the tortoises  \label{calc_area_tort}}

For the coordinates $x (\phi), y (\phi)$ of the boundary of the sets, see \ref{x_phi} and \ref{y_phi}. Now the boundary is rotated by $- \phi$. In the case $D_0$ the boundary is a circle and the
point $(x (\phi), y (\phi))$ is rotated to $(1, 0)$. The tangent on the point $(x (\phi), y (\phi))$ rotated
by $- \phi$ is vertical. So the formulae belonging to figure \ref{gseg2} can be applied.
\begin{equation}
   \begin{split}
   & x_\phi (\phi) = + \cos (\phi) \, x (\phi) + \sin (\phi) \, y (\phi) \quad
   \bar{x}_\phi (\phi) = x_\phi (\phi) - 1 \\
   & y_\phi (\phi) = - \sin (\phi) \, x (\phi) + \cos (\phi) \, y (\phi)
   \label{xy_rot}
   \end{split}
\end{equation}
$\bar{x}_\phi (\phi), y_\phi (\phi)$ have the property:
\begin{equation}
   \bar{x}_\phi (\phi + \pi) = \bar{x}_\phi (\phi)  \qquad y_\phi (\phi + \pi) = x_\phi (\phi) 
   \label{xy_rot_prop}
\end{equation}

With $\psi = \pi / 3$ we get, all quantities are linear in $\epsilon$:
\begin{equation}
   \begin{split}
   & d_{x,k} = \bar{x}_\phi (2 \, k \, \psi) + \bar{x}_\phi ((2 \, k + 1) \, \psi) \\
   & d_{y,k} =       y_\phi (2 \, k \, \psi) +       y_\phi ((2 \, k + 1) \, \psi) \\
   & r_{lu,k} = - q (2 \, k \, \psi)_+  \qquad \qquad r_{ll,k} = - q (2 \, k \, \psi)_- \\
   & r_{ru,k} = - q ((2 \, k + 1) \, \psi)_+  \quad r_{rl,k} = - q ((2 \, k + 1) \, \psi)_-
   \label{xy_rot_4r}
   \end{split}
\end{equation}
$q (\phi)$ is discontinuous, the subscripts $\pm$ above denote the left and right values at a discontinuity. The 3 subscripts $k = 0, 1, 2$ correspond to the configuration belonging to the edges $E_k$ in figure \ref{lattice_orientations}.  \\

For $\hat{A_2} (\dots)$ i.e. 2-parameter minimizing, see \ref{A2hat_d_4r}, we get for the area of the $6$ segments, cutted from $D_\epsilon$:
\begin{equation}
   \sum\limits_{k=0}^{2} \hat{A_2} (d_{x,k}, r_{lu,k}, r_{ll,k}, r_{ru,k}, r_{rl,k}) 
   \label{area_6_cutted}
\end{equation}
Because of the property \ref{q_phi_prop} for $q (\phi)$ and the properties \ref{xy_rot_prop} for $\bar{x}_\phi (\phi), y_\phi (\phi)$, the linear term in $\epsilon$ is $0$. \\

As final result the area of tortoise is, now depending on $\epsilon$:
\begin{equation}
   A_T = \underbrace{\pi - 0.0104747054 \, \epsilon^2}_{\text{area of} \, D_\epsilon, \, \text{see \ref{xy_area}}}
   \quad - \quad \underbrace{(6 \, A_1 (0,0) - 0.0118673317 \, \epsilon^2 + O(\epsilon^3))}_{\sum A_2 (\dots)} 
   \nonumber
\end{equation}
\begin{equation}
   = \underbrace{\pi - 6 \, A_1 (0,0)}_{\text{Croft's tortoise}} \quad + \quad 0.0013926262 \, \epsilon^2
   \quad + \quad O(\epsilon^3)
   \label{area_tortoise}
\end{equation}

Because of a $+$ sign for the coefficient of $\epsilon^2$, the density of the 2-avoiding set $S_\epsilon$ has a local minima for $\epsilon = 0$. The density at the minima is the density of Croft's tortoises. \qedsymbol  \\

In this case for the given $D_\epsilon$, 2-parameter minimizing seems to be essential.  This because for  1-parameter minimizing using  $A_2 ()$, see \ref{A2_d_4r}, it is not clear, if Croft's tortoise has a local minima/maxima for the area (numerical accuracy). In the area now the trailing term is very small $A_T = \dots + 2.04 10^{-15}  \, \epsilon^2 + O(\epsilon^3)$. \\

\section{Conclusion}

In this article I constructed a 1-parameter family $S_\epsilon$ of 2-avoiding sets. For small $\epsilon$ the sets have  a density $\ge$ the density of Croft's tortoises, which are unit distance avoiding. The possibilities for errors in this article are numerous. Therefore it would be interesting to know
the exact density (within the numerically possible accuracy) for small $\epsilon$. If this density
would be greater than Croft's density, this would be an evidence but not more, for this paper to be correct.
It would be interesting too, for which $\epsilon$ the density has its maximum and its numerical value.

\newpage
\noindent \textbf{\large Appendices}

\appendix

\section{Croft's construction of a planar, unit distance \\ avoiding set using tortoises}

Croft places copies of open discs with diameter $1.0$ on a hexagonal lattice with lattice constant $L = 2.0$, i.e. the 6 shortest basis vectors in this lattice have length $L$.
This set is invariant under the hexagonal group and is a unit distance avoiding set with density
 $\pi \, / \, (8 \, \sqrt{3})$. Decreasing $L$ the set is no longer unit distance avoiding. See the following figure: we cut the 2 blue disc segments in a vertical strip of width $1.0$, this for each of the $6$ directions. The remaining set is again unit distance avoiding. To get the maximal density we have to choose $L_C = 1.96553$. The cutted segment has $\phi_C = 15.08686^{\circ}$ as half segment angle, so $L_C = 1 + \cos (\phi_C)$. See also figure 2. in \cite{CmCsMaVaZs}.  

\begin{center}
  \includegraphics[width=0.92\textwidth]{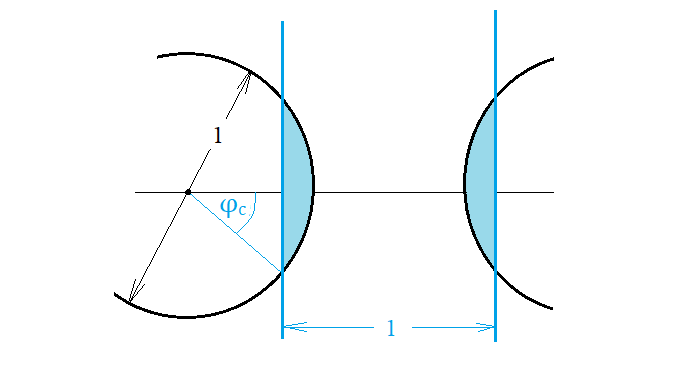}
  \captionof{figure}{Croft's disc segment in blue.}
  \label{croft}
\end{center}

In the rest of this paper Croft's segment is always used for diameter $2$. So here the numerial values
for quantities, describing this segment.  $\phi_C = 0.263315538964831$ remains the half segment angle, the horizontal width of the segment is $w_C = 1 - \cos (\phi_C) = 0.034467692551095$ and the area is $A_C = \phi_C - \cos (\phi_C) \, \sin (\phi_C) = 0.012003664907850$.

\section{The area of a general disc segment as power series \label{sect_gseg}}

Croft's construction uses discs of \emph{diameter} $1$ at the lattice points and so the periodic set is a unit distance or 1-avoiding set. Here we use discs with \emph{radius} $1$. So we construct 2-avoiding sets. But this has no influence on the density because the lattice constant is doubled too.

\begin{center}
  \includegraphics[width=0.50\textwidth]{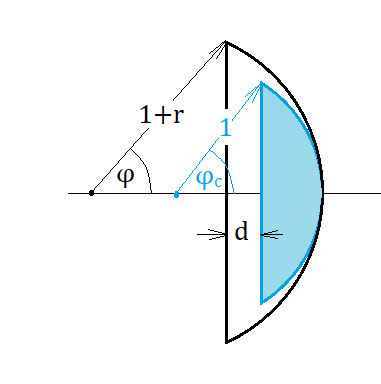}
  \captionof{figure}{The general disc segment in black. For $d, r = 0$ we get
                     Croft's disc segment in blue.}
  \label{gseg}
\end{center}

The area $A_1$ of the black disc segment is given by:
\begin{equation}
   R = 1 + r \qquad D = 1 - \cos (\phi_C) + d \qquad \phi = \arccos ((R - D) / R))
   \nonumber
\end{equation}
\begin{equation}
   A_1 (d, r) = R^2 \, \phi - (R - D) \, R \sin (\phi)
   \label{def_D_3}
\end{equation}

Now we give the area $A_1 (d, r)$ as power series in $d, r$. The first and second derivatives in the following are taken at $d, r = 0$.
\begin{equation}
   \begin{split}
   & B = \; \; \frac{dA_1}{dd}   \Big\vert_{d,r=0} = 2 \, \sin (\phi_C)            = + 0.5205664 \\
   & C = \; \; \frac{dA_1}{dr}   \Big\vert_{d,r=0} = 2 \, (\phi_C - \sin (\phi_C)) = + 0.0060646 \\ 
   & D = \frac{d^2A_1}{dd^2}  \Big\vert_{d,r=0}    = 2 \, \cot (\phi_C)            = + 7.4190894 \\
   & E = \frac{d^2A_1}{dd \, dr} \Big\vert_{d,r=0} = 2 \, \frac{1 - \cos (\phi_C)}{\sin (\phi_C)}
                                                                                   = + 0.2648475 \\
   & F = \frac{d^2A_1}{dr^2}  \Big\vert_{d,r=0}    = 2 \, \left( \phi_C - 2 \, \frac{1 - \cos (\phi_C)}
                                                         {\sin (\phi_C)} \right) = - 0.0030640 
   \label{FiSecDer}
   \end{split}
\end{equation}

The power series for $A_1 (d, r)$ is now up to terms of second degree:
\begin{equation}
   A_1 (d, r) = A_1 (0, 0) + B \, d + C \, r + D / 2 \, d^2 + E \, d \, r + F / 2 \, r^2 + \dots
   \label{PowSerA1}
\end{equation}

\begin{center}
  \includegraphics[width=0.75\textwidth]{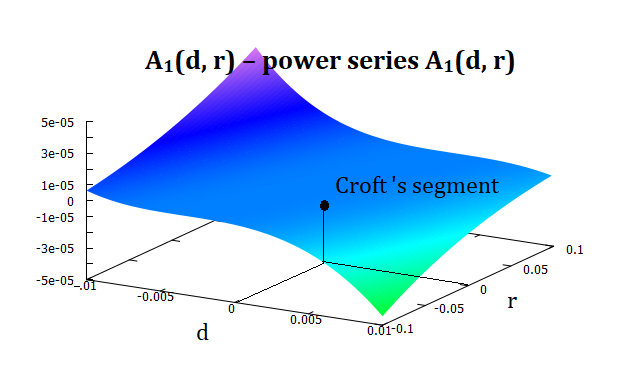}
  \captionof{figure}{The difference between the exact value of $A_1 (d, r)$ and the approximation
                     by the power series \ref{PowSerA1} for $-0.01 \le d \le +0.01, \;
                     -0.1 \le r \le +0.1$. Croft's segment has a width of $1 - \cos (\phi_C) = 0.0344$,
                     so the negative range of $d$ covers approx. $30\%$ of the segment. }
  \label{DifA1}
\end{center}

\section{The area of two general, opposite disc segments as power series \label{2_opp_segs}}

\begin{center}
  \includegraphics[width=0.82\textwidth]{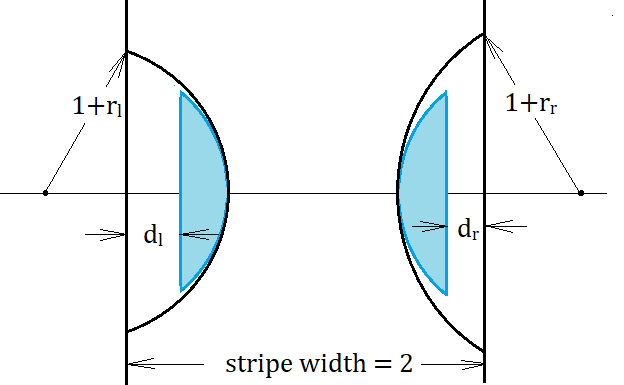}
  \captionof{figure}{The general two disc segments in black. For $d_l, d_r, r_l, r_r = 0$ we get
                     Croft's opposite disc segments in blue, see figure \ref{croft}.}
  \label{gseg2}
\end{center}

The vertical stripe of width $2$ in the figure above can be shifted in horizontal direction by $s$
without changing the property of the set being 2-avoiding. We choose the $s$ for which the sum of
the areas of the two disc segments $\bar{A_2} (d_l, d_r, r_l, r_r, \; s)$ is a minimum. For the coefficients
$B, \dots, F$ see \ref{PowSerA1}.
\begin{equation}
   \begin{split}
   \bar{A_2} & (d_l, d_r, r_l, r_r, \; s) = \\
                              & A_1 (0, 0) + B \, (d_l - s) + C \, r_l + D / 2 \, (d_l - s)^2 +
                                E \, (d_l - s) \, r_l + F / 2 \, r_l^2 + \\
                              & A_1 (0, 0) + B \, (d_r + s) + C \, r_r + D / 2 \, (d_r + s)^2 +
                                E \, (d_r + s) \, r_r + F / 2 \, r_r^2 + \dots
   \label{Min_A2bar_s}
   \end{split}
\end{equation}

This $s_{min}$ is determined by the equation $\frac{d\bar{A_2} (\dots, \; s)}{ds} = 0$. Inserting $s_{min}$ in \ref{Min_A2bar_s} above we get the following result. Now $\bar{A_2}$ does not depend on $d_l, d_r$ but only on the sum $d = d_l + d_r$:
\begin{equation}
   \begin{split}
   \bar{A_2} & (d, r_l, r_r) = \; 2 \, A_1 (0, 0) + B \, d + C \, (r_l + r_r) + \frac{D}{4} \, d^2  \\
                       & + \frac{2 \, D \, F - E^2}{4 \, D} \, (r_l^2 + r_r^2)
                       + \frac{E}{2} \, d \, (r_l + r_r) + \frac{E^2}{2 \, D} \,  r_l \, r_r + \dots 
   \label{A2_d_rlrr}
   \end{split}
\end{equation}
\begin{equation}
   \begin{split}
   \bar{A_2} (d, r_l, r_r) = \; & 2 \, A_1 (0, 0) + 0.5205664 \, d + 0.0060646 \, (r_l + r_r) + 1.8547723 \, d^2  \\
                       & - 0.0038956 \, (r_l^2 + r_r^2)
                       + 0.1324237 \, d \, (r_l + r_r) + 0.0047272 \,  r_l \, r_r + \dots 
   \label{A2_d_rlrr_num}
   \end{split}
\end{equation}

The formula above is for the case that we have only one radius at one side in \ref{gseg2}. The radii of curvature for the sets of constant diameter used in this article are discontinuous at a finite number of points. Therefore also a formula for $2$ radii $r_{lu}, r_{ll}$ on the left and $2$ radii $r_{ru}, r_{rl}$ on the right is necessary. The second subscript means (u)pper and (l)ower. \\

The upper and lower disc segments left and right contribute with the factor $1/2$ to the area to minimize.
With the minimizer of above \ref{Min_A2bar_s} this area is:
\begin{equation}
   \begin{split}
   A_2 (d, r_{lu}, r_{ll}, r_{ru}, r_{rl}, \enspace s)
       = ( \; & A_1 (d / 2 {\color{red}  \bf \large \, +}  \,  s, r_{lu} )
                           + A_1 (d / 2 {\color{blue}  \bf \large \, -}  \, s, r_{ru} ) \\
            + & A_1 (d / 2 {\color{red}  \bf \large \, +}  \,  s, r_{ll} ) \,
                           + A_1 (d / 2 {\color{blue}  \bf \large \, -}  \, s, r_{rl} ) \; ) \, / \, 2
   \label{Min_A2_s}
   \end{split}
\end{equation}

After minimizing, using $r_s = r_{lu} + r_{ll} + r_{ru} + r_{rl}, r_{s2} = r_{lu}^2 + r_{ll}^2 + r_{ru}^2 + r_{rl}^2,$ $r_l = r_{lu} + r_{ll} - r_{ru} - r_{rl}$ as short forms:
\begin{equation}
   \begin{split}
   A_2 (d, r_{lu}, r_{ll}, r_{ru}, r_{rl}) = & \; 2 \, A_1 (0, 0)
         \; {\color{red} + B \, d + \frac{C}{2} \, r_s} \\
         &  + \frac{D}{4} \, d^2 + \frac{E}{4} \, d \, r_s - \frac{E^2}{16 \, D} \,  r_l^2 + \frac{F}{4} \, r_{s2} + \dots 
   \label{A2_d_4r}
   \end{split}
\end{equation}
When summing over the $3$ orientations of opposite $D_\epsilon$, see figure \ref{lattice_orientations}, the red, linear terms above cancel and the result is the desired quadratic form plus twice the area of Croft's segment.
\begin{equation}
   \begin{split}
   A_2 & (d, r_{lu}, r_{ll}, r_{ru}, r_{rl}) = \; 2 \, A_1 (0, 0) + 0.5205664 \, d
           + 0.0030323 \, r_s   \\
         & + 1.8547723 \, d^2 + 0.0662118 \, d \, r_s - 0.0005909 \, r_l^2 - 7.6602006 \, r_{s2} + \dots 
   \label{A2_d_4r_num}
   \end{split}
\end{equation}

\begin{remark}  \label{rem_min2}
When minimizing the area of two opposite caps, only the horizontal position of the stripe is varied.
This is enough, when the configuration of the disks has a up-down reflection symmetry as in figure \ref{gseg2}.
Varying the inclination of the stripe too, would lead to increasing areas. But in the case \ref{Min_A2_s}
with different upper and lower radii, this up-down symmetry does not exist. Varying the inclination too,
could result in an even lower area of the segments to cut. See appendix \ref{2par_area_opp2}. 
\end{remark}

\section{The 3-color vertex coloring of a hexagonal lattice \label{3col_hex_lattice}}

\tikzstyle{bullet_red}   = [circle, color = red,   draw, fill, inner sep = 0 pt, minimum size = 2.0 mm]
\tikzstyle{bullet_green} = [circle, color = green, draw, fill, inner sep = 0 pt, minimum size = 2.0 mm]
\tikzstyle{bullet_blue}  = [circle, color = blue,  draw, fill, inner sep = 0 pt, minimum size = 2.0 mm]

\begin{figure}[htp]
\centerline{
\begin{tikzpicture}[scale = 1]
\def \s{2.3}                     
\def \h{0.86602}                 
\def \ared  {30}                 
\def \agreen{\ared + 30}         
\def \ablue {\ared + 60}         

\foreach \i in { 0, ..., 3 }  
   \draw[color = black] ( 0.5 * \i * \s, \i * \h * \s ) --  ( 3.0 * \s + 0.5 * \i * \s, \i * \h * \s ); 
\foreach \i in { 0, ..., 3 }  
   \draw[color = black] ( 1.0 * \i * \s,  0 * \h * \s ) --  ( 1.5 * \s + 1.0 * \i * \s,  3 * \h * \s ); 
\foreach \i in { 1, ..., 3 }  
   \draw[color = black] ( 0.5 * \i * \s           , 1 * \i * \h * \s ) --  ( 1.0 * \i * \s, 0  * \h * \s ); 
\foreach \i in { 1, ..., 2 }  
   \draw[color = black] ( 1.5 * \s + 1.0 * \i * \s, 3 *      \h * \s ) --  ( 3.0 * \s + 0.5 * \i * \s, \i * \h * \s ); 

\draw[color = red, thick] ( 0.0 * \s, 0 * \h * \s ) -- ( 3.0 * \s, 0 * \h * \s ) --
                          ( 4.5 * \s, 1 * \h * \s ) -- ( 1.5 * \s, 1 * \h * \s ) -- cycle;

\draw[color = blue, thick] ( 0.5 * \s, 1 * \h * \s ) -- ( 2.0 * \s, 0 * \h * \s ) --
                           ( 3.5 * \s, 1 * \h * \s ) -- ( 2.0 * \s, 2 * \h * \s ) -- cycle;

\node[draw = none, below left]()   at ( 0.0 * \s, 0.0 * \s )  {\Large $0$};


\node at ( 0.0 * \s, 0 * \h * \s ) [bullet_red]   {};
\node at ( 1.0 * \s, 0 * \h * \s ) [bullet_green] {};
\node at ( 2.0 * \s, 0 * \h * \s ) [bullet_blue]  {};
\node at ( 3.0 * \s, 0 * \h * \s ) [bullet_red]   {};

\node at ( 0.5 * \s, 1 * \h * \s ) [bullet_blue]  {};
\node at ( 1.5 * \s, 1 * \h * \s ) [bullet_red]   {};
\node at ( 2.5 * \s, 1 * \h * \s ) [bullet_green] {};
\node at ( 3.5 * \s, 1 * \h * \s ) [bullet_blue]  {};

\node at ( 1.0 * \s, 2 * \h * \s ) [bullet_green] {};
\node at ( 2.0 * \s, 2 * \h * \s ) [bullet_blue]  {};
\node at ( 3.0 * \s, 2 * \h * \s ) [bullet_red]   {};
\node at ( 4.0 * \s, 2 * \h * \s ) [bullet_green] {};

\node at ( 1.5 * \s, 3 * \h * \s ) [bullet_red]   {};
\node at ( 2.5 * \s, 3 * \h * \s ) [bullet_green] {};
\node at ( 3.5 * \s, 3 * \h * \s ) [bullet_blue]  {};
\node at ( 4.5 * \s, 3 * \h * \s ) [bullet_red]   {};

\node at ( 4.5 * \s, 1 * \h * \s ) [bullet_red] {};

\node[draw = none, color = black, below      ]()  at ( 0.5  * \s, 0   * \h * \s - 0.05 * \s )   {\large $\omega_1$};
\node[draw = none, color = black, above left ]()  at ( 0.25 * \s, 0.5 * \h * \s )   {\large $\omega_2$};

\end{tikzpicture}
}
\caption{\label{3-Coloring}
The unique vertex coloring of the hexagonal lattice (with respect to permutations of the colors).
Let $\omega_1$ be the horizontal vector (red origin to green) with the length of a triangle side, $\omega_2$ the vector in the direction $2 \, \pi / 3$ (red origin to blue) also with the length of a triangle side.}
\end{figure}
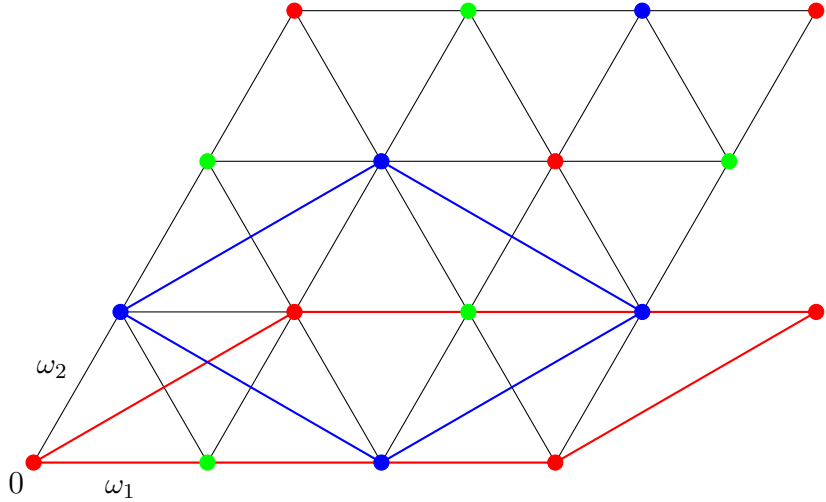

The symmetry group $C_f$ of isometric transformations of $\R^2$ fixing a color is generated by $\omega_1,
\omega_1 +  \omega_2$, a rotation around the origin at 0 by an angle of $2 \pi / 3$ of order $3$ and a reflection on a horizontal line through the origin. This is the wallpaper group $p31m$.

\begin{equation}
   C_f \simeq (\Z \times \Z) \; \rtimes  \underbrace{(\Z_3 \rtimes \Z_2)}_{S_3}
 \label{ColorFix} \end{equation}
where $\rtimes$ denotes the semidirect product. 
One fundamental parallelogram of the red lattice can be seen in red in figure \ref{3-Coloring}.
Another normalized (closer to a rectangle) fundamental parallelogram is the rhombus of blue lattice vertices. \\

The symmetry group $C_p$ with additional transformations permuting the 3 colors is generated by an
 additional shift by $\omega_1$ (a cyclic permutation on the 3 colors) and a rotation $\pi / 6$
 (a transposition on the 3 colors). : denotes a non-split extension.
\begin{equation}
   C_p \simeq C_f : \underbrace{(\Z_3 \times \Z_2)}_{C_6}
 \label{ColorPerm} \end{equation}

\tikzstyle{bullet_red}    = [circle, color = black!50,  draw, fill, inner sep = 0 pt, minimum size = 1.5 mm]
\tikzstyle{bullet_green}  = [circle, color = black!50,  draw, fill, inner sep = 0 pt, minimum size = 1.5 mm]
\tikzstyle{bullet_blue}   = [circle, color = black!50,  draw, fill, inner sep = 0 pt, minimum size = 1.5 mm]

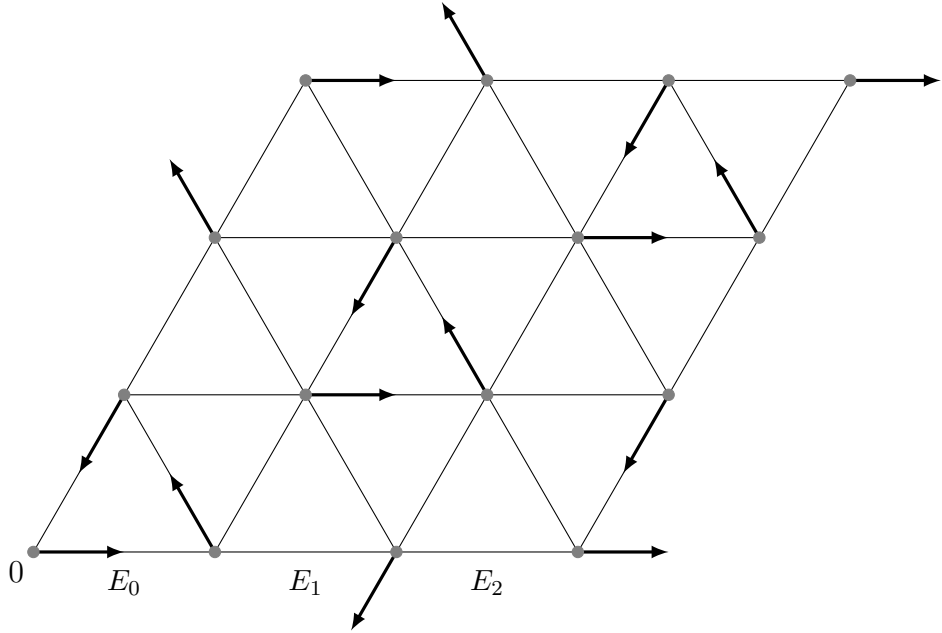
\begin{figure}[htp]
\centerline{
\begin{tikzpicture}[scale = 1]
\def \s{2.4}                     
\def \h{0.86602}                 
\def \rpd{1.2}                   
\def \ared  {30}                 
\def \agreen{\ared + 30}         
\def \ablue {\ared + 60}         

\foreach \i in { 0, ..., 3 }  
   \draw[color = black] ( 0.5 * \i * \s, \i * \h * \s ) --  ( 3.0 * \s + 0.5 * \i * \s, \i * \h * \s ); 
\foreach \i in { 0, ..., 3 }  
   \draw[color = black] ( 1.0 * \i * \s,  0 * \h * \s ) --  ( 1.5 * \s + 1.0 * \i * \s,  3 * \h * \s ); 
\foreach \i in { 1, ..., 3 }  
   \draw[color = black] ( 0.5 * \i * \s           , 1 * \i * \h * \s ) --  ( 1.0 * \i * \s, 0  * \h * \s ); 
\foreach \i in { 1, ..., 2 }  
   \draw[color = black] ( 1.5 * \s + 1.0 * \i * \s, 3 *      \h * \s ) --  ( 3.0 * \s + 0.5 * \i * \s, \i * \h * \s );

\node[draw = none, below left]()   at ( 0.0 * \s, 0.0 * \s )  {\Large $0$};

\def \arl{0.5 * \s}   
\def \arx{0.}         
\def \ary{0}          
\def \ardx{0}         
\def \ardy{0}         

\begin{scope}[>= latex]

   \def \ardx{1.00000 * \arl} \def \ardy{0.00000 * \arl}    

   \def \arx{0.0 * \s}  \def \ary{0 * \h * \s}
   \draw[->, color = black, very thick]  (\arx, \ary) to (\arx + \ardx, \ary + \ardy);
   \def \arx{1.5 * \s}  \def \ary{1 * \h * \s}
   \draw[->, color = black, very thick]  (\arx, \ary) to (\arx + \ardx, \ary + \ardy);
   \def \arx{3.0 * \s}  \def \ary{2 * \h * \s}
   \draw[->, color = black, very thick]  (\arx, \ary) to (\arx + \ardx, \ary + \ardy);
   \def \arx{4.5 * \s}  \def \ary{3 * \h * \s}
   \draw[->, color = black, very thick]  (\arx, \ary) to (\arx + \ardx, \ary + \ardy);
   \def \arx{3.0 * \s}  \def \ary{0 * \h * \s}
   \draw[->, color = black, very thick]  (\arx, \ary) to (\arx + \ardx, \ary + \ardy);
   \def \arx{1.5 * \s}  \def \ary{3 * \h * \s}
   \draw[->, color = black, very thick]  (\arx, \ary) to (\arx + \ardx, \ary + \ardy);

   \def \ardx{- 0.50000 * \arl} \def \ardy{0.86602 * \arl}    

   \def \arx{1.0 * \s}  \def \ary{0 * \h * \s}
   \draw[->, color = black, very thick]  (\arx, \ary) to (\arx + \ardx, \ary + \ardy);
   \def \arx{2.5 * \s}  \def \ary{1 * \h * \s}
   \draw[->, color = black, very thick]  (\arx, \ary) to (\arx + \ardx, \ary + \ardy);
   \def \arx{4.0 * \s}  \def \ary{2 * \h * \s}
   \draw[->, color = black, very thick]  (\arx, \ary) to (\arx + \ardx, \ary + \ardy);
   \def \arx{1.0 * \s}  \def \ary{2 * \h * \s}
   \draw[->, color = black, very thick]  (\arx, \ary) to (\arx + \ardx, \ary + \ardy);
   \def \arx{2.5 * \s}  \def \ary{3 * \h * \s}
   \draw[->, color = black, very thick]  (\arx, \ary) to (\arx + \ardx, \ary + \ardy);

   \def \ardx{- 0.50000 * \arl} \def \ardy{- 0.86602 * \arl}    

   \def \arx{2.0 * \s}  \def \ary{0 * \h * \s}
   \draw[->, color = black, very thick]  (\arx, \ary) to (\arx + \ardx, \ary + \ardy);
   \def \arx{3.5 * \s}  \def \ary{1 * \h * \s}
   \draw[->, color = black, very thick]  (\arx, \ary) to (\arx + \ardx, \ary + \ardy);
   \def \arx{0.5 * \s}  \def \ary{1 * \h * \s}
   \draw[->, color = black, very thick]  (\arx, \ary) to (\arx + \ardx, \ary + \ardy);
   \def \arx{2.0 * \s}  \def \ary{2 * \h * \s}
   \draw[->, color = black, very thick]  (\arx, \ary) to (\arx + \ardx, \ary + \ardy);
   \def \arx{3.5 * \s}  \def \ary{3 * \h * \s}
   \draw[->, color = black, very thick]  (\arx, \ary) to (\arx + \ardx, \ary + \ardy);

\end{scope}

\node at ( 0.0 * \s, 0 * \h * \s ) [bullet_red]   {};
\node at ( 1.0 * \s, 0 * \h * \s ) [bullet_green] {};
\node at ( 2.0 * \s, 0 * \h * \s ) [bullet_blue]  {};
\node at ( 3.0 * \s, 0 * \h * \s ) [bullet_red]   {};

\node at ( 0.5 * \s, 1 * \h * \s ) [bullet_blue]  {};
\node at ( 1.5 * \s, 1 * \h * \s ) [bullet_red]   {};
\node at ( 2.5 * \s, 1 * \h * \s ) [bullet_green] {};
\node at ( 3.5 * \s, 1 * \h * \s ) [bullet_blue]  {};

\node at ( 1.0 * \s, 2 * \h * \s ) [bullet_green] {};
\node at ( 2.0 * \s, 2 * \h * \s ) [bullet_blue]  {};
\node at ( 3.0 * \s, 2 * \h * \s ) [bullet_red]   {};
\node at ( 4.0 * \s, 2 * \h * \s ) [bullet_green] {};

\node at ( 1.5 * \s, 3 * \h * \s ) [bullet_red]   {};
\node at ( 2.5 * \s, 3 * \h * \s ) [bullet_green] {};
\node at ( 3.5 * \s, 3 * \h * \s ) [bullet_blue]  {};
\node at ( 4.5 * \s, 3 * \h * \s ) [bullet_red]   {};

\node[draw = none, color = black, below      ]()  at ( 0.5  * \s, 0   * \h * \s - 0.05 * \s )   {\large $E_0$};
\node[draw = none, color = black, below      ]()  at ( 1.5  * \s, 0   * \h * \s - 0.05 * \s )   {\large $E_1$};
\node[draw = none, color = black, below      ]()  at ( 2.5  * \s, 0   * \h * \s - 0.05 * \s )   {\large $E_2$};


\end{tikzpicture}
}
\caption{\label{lattice_orientations}
The hexagonal lattice decorated with vectors, showing the 3 orientations of the copies of the constant diameter sets $D_\epsilon$. The graph has 3 different types of cuts on the edges. The group of isometries, acting on the lattice and conserving the orientations, acts in 3 different orbits on the edges $E_0, E_1, E_2$, representing the orbits. These 3 edeges also represent the 3 types of cuts with stripes. The group acts transitive on the vertices and in 6 orbits on the triangles.}
\end{figure}

\newpage

\section{Copy \& Paste into a Computer Algebra System \label{copa_1_par}}

The piecewise linear functions are defined on these intervals, the interval ends are not equidistant. The values have to be multiplied by $\pi$: \\
{ \ttfamily \footnotesize
0,4/45,1/6, \; 11/45,1/3,19/45, \; 1/2,26/45,2/3, \; 34/45,5/6,41/45, \; 1,49/45,7/6, \\
56/45,4/3,64/45, \;  3/2,71/45,5/3, \;  79/45,11/6,86/45, \; 2} \\

\noindent
$24$ values for $q (\phi)$ with the property $q_{i+12} = - q_i$, see \ref{q_phi}: \\
{ \ttfamily \footnotesize
-0.977901957024321,+0.724209871347166,-0.733569955967565,+1.000000000000000, \\
-0.922743876968233,+0.488920844394468,-0.126049416295258,+0.044264839209546, \\
+0.004202409557006,-0.101935908145429,+0.472228160625608,-0.899297801582359, \\

\noindent
+0.977901957024321,-0.724209871347166,+0.733569955967565,-1.000000000000000, \\
+0.922743876968233,-0.488920844394468,+0.126049416295258,-0.044264839209546, \\
-0.004202409557006,+0.101935908145429,-0.472228160625608,+0.899297801582359  \\
} \\

\noindent
$24$ values for the piecewise part of $x (\phi)$ with the property $x_{m,i+12} = + x_{m,i}$, see \ref{x_phi}: \\
{ \ttfamily \footnotesize
-0.977901957024321,+0.658272945792727,-0.604201417786322,+0.642824448212470, \\
-0.318547490271646,+0.022965115071595,+0.022965115071595,-0.018237632467773, \\
+0.001793582358496,+0.078143098622847,-0.419097570873107,+0.899297801582358, \\

\noindent
-0.977901957024321,+0.658272945792727,-0.604201417786322,+0.642824448212470, \\
-0.318547490271646,+0.022965115071595,+0.022965115071595,-0.018237632467773, \\
+0.001793582358496,+0.078143098622847,-0.419097570873107,+0.899297801582358  \\
} \\

\noindent
$24$ values for the piecewise part of $y (\phi)$ with the property $y_{m,i+12} = + y_{m,i} $, see \ref{y_phi}: \\
{ \ttfamily \footnotesize
+0.000000000000000,+0.469165603677154,-0.259724309980211,+0.944514570708893, \\
-0.720630471716577,+0.649101774356247,+0.034131513666520,+0.199386707906710, \\
+0.164691626090284,+0.090961755271850,+0.378043789657369,+0.000000000000000, \\

\noindent
+0.000000000000000,+0.469165603677154,-0.259724309980211,+0.944514570708893, \\
-0.720630471716577,+0.649101774356247,+0.034131513666520,+0.199386707906710, \\
+0.164691626090284,+0.090961755271850,+0.378043789657369,+0.000000000000000  \\
} \\

Before rotating by $2 \pi / 3, 4 \pi / 3$ the set has to be shifted in x/y-direction by:
\begin{equation}
   {\Delta}x = \text{\ttfamily \footnotesize - 0.001383301426275} \; \epsilon  \qquad
   {\Delta}y = \text{\ttfamily \footnotesize - 0.158574235421304} \; \epsilon
   \label{shift}
\end{equation}

\noindent
The lattice constant of the hexagonal lattice is $L = \,$ { \ttfamily \footnotesize 3.93106461489781} \\
The area of the set $D_\epsilon$ is $\pi \; - $ { \ttfamily \footnotesize 0.010474705472633 } $\epsilon^2$.

\begin{center}
  \includegraphics[width=0.50\textwidth]{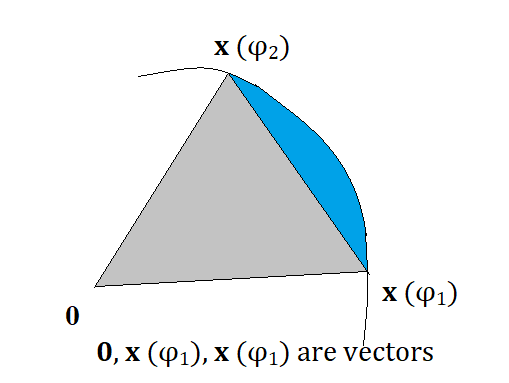}
  \captionof{figure}{Calculating the area of a segment. It is the difference between the area of the sector
                     (blue and grey) and the grey triangle.}
  \label{calcseg}
\end{center}

It is simple to calculate the area of the blue segment in the figure above with corners $\mathbf{x} (\phi_1), \mathbf{x} (\phi_2)$ (bounded by a line and an arc of the the curve through these points)
\begin{equation}
   \underbrace{\frac{1}{2} \int\limits_{\phi_1}^{\phi_2} \left( x (\phi) \, \frac{dy (\phi)}{d\phi} 
                - y (\phi) \, \frac{dx (\phi)}{d\phi}  \right) d\phi}
     _{\text{area of the sector} \; \mathbf{0} - \mathbf{x} (\phi_1) - \mathbf{x} (\phi_2) }
   \quad - \quad
   \underbrace{\frac{1}{2} \left( x (\phi_1) \, y (\phi_2) - x (\phi_2) \, y (\phi_1) \right) }
     _{\text{area of the triangle} \; \mathbf{0} - \mathbf{x} (\phi_1) - \mathbf{x} (\phi_2) }
   \label{area_segment calc}
\end{equation}

\section{2-parameter minimizing, the area of a more general disc segment as power series  \label{3par_segm}}

This is an extension of the disc segment \ref{gseg}, defined in appendix \ref{sect_gseg} by an additional 
parameter $\delta$. This parameter is needed in the case, when the inclination of the cutted stripe is varied too.

\begin{center}
  \includegraphics[width=0.50\textwidth]{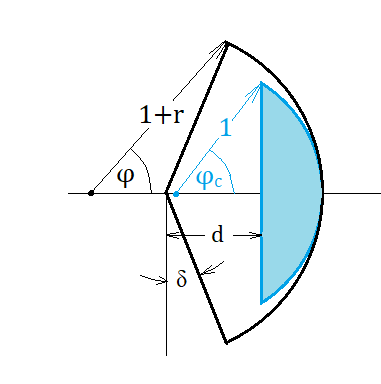}
  \captionof{figure}{The general, 3-parameter disc segment in black. For $d, r, \delta = 0$ we get
                     Croft's disc segment in blue.}
  \label{gseg2par}
\end{center}

The area $\bar{A_1}$ of the black disc segment is given by:
\begin{equation}
   R = 1 + r \qquad D = 1 - \cos (\phi_C) + d \qquad \cos (\phi) = (R - D + R \sin (\phi) \tan (\delta)) / R
   \nonumber
\end{equation}
\begin{equation}
   \phi = \arccos ((R - D) / R \, \cos (\delta)) - \delta
   \label{def_phi}
\end{equation}
\begin{equation}
   \bar{A_1} (d, r, \delta) = R^2 \, \phi - (R - D) \, R \sin (\phi)
   \label{def_D_3}
\end{equation}

Now we give the area $\bar{A_1} (d, r, \delta)$ as power series in $d, r, \delta$. The first and second derivatives in the following are taken at $d, r, \delta = 0$.
\begin{equation}
   \begin{split}
   & H = \; \, \frac{d\bar{A_1}}  {d\delta}       \Big\vert_{d,r,\delta=0} = - \sin (\phi_C)^2
                                                                                     = - 0.0677473 \\
   & J = \; \frac{d^2\bar{A_1}}{dd \, d\delta} \Big\vert_{d,r,\delta=0} = - 2 \, \cos (\phi_C)
                                                                                     = - 1.9310646 \\
   & K = \frac{d^2\bar{A_1}}{dr \, d\delta}    \Big\vert_{d,r,\delta=0} = 2 \, (\cos (\phi_C) - 1)
                                                                                     = - 0.0689353 \\
   & L = \, \frac{d^2\bar{A_1}}{d\delta^2}        \Big\vert_{d,r,\delta=0}  = 2 \, \cos (\phi_C) \, \sin (\phi_C)
                                                                                     = + 0.5026237 
   \label{A1_hat_SecDer}
   \end{split}
\end{equation}

The power series for $\bar{A_1} (d, r, \delta)$ is now up to terms of second degree, $A_1 (d, r)$ was already given in appendix \ref{sect_gseg}:
\begin{equation}
   \bar{A_1}  (d, r, \delta) = A_1 (d, r) \quad + H \, \delta + J \, d \, \delta
                                                + K \, r \, \delta +  L / 2 \, \delta^2 + \dots
   \label{PowSerA1_hat}
\end{equation}

\begin{center}
  \includegraphics[width=0.80\textwidth]{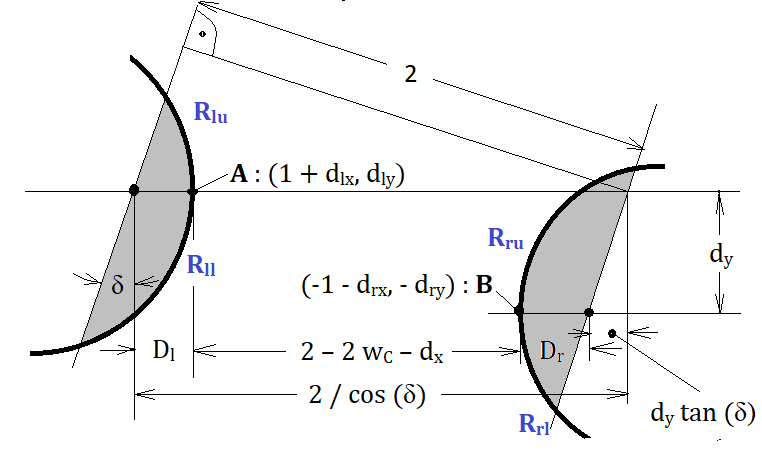}
  \captionof{figure}{The most general two disc segments in black. The 4 radii in blue can be different.
                     The displacements of the points $A$ and $B$ are given relative to the corresponding
                     lattice points at the left and right. The lattice constant is
                     $4 - 2 \, w_C$. $w_C = 1 - \cos (\phi_C)$ is the horizontal width of
                     Croft's disc segment. $d_x = d_{lx} + d_{rx} , d_y = d_{ly} + d_{ry}$.}
  \label{gseg2_exact}
\end{center}

2-parameter minmizing with $\delta$ needs an exact $\bar{d}_x$ including $\delta$, not only the $d_{x(,k)}$ given in \ref{xy_rot_4r}. This exact $\bar{d}_x$ can be extracted from figure \ref{gseg2_exact} above. \\
\begin{equation}
   \begin{split}
   & \bar{d}_x = d_x + 2 \, (1 / \cos (\delta) - 1) - \tan (\delta) \, d_y \\  
   & \bar{d}_x = d_x + \delta^2 - \delta \, d_y \quad + \dots \qquad  \text{terms up to second order}
   \label{d_exact}
   \end{split}
\end{equation}

The upper and lower disc segments left and right contribute with the factor $1/2$ to the area to minimize.
This area is:
\begin{equation}
   \begin{split}
   \hat{A_2} (\bar{d}_x, r_{lu}, r_{ll}, r_{ru}, & r_{rl}, \enspace s, \delta) \\
       = ( \; & \bar{A_1} (\bar{d}_x / 2 {\color{red}  \bf \large \, +}  \,  s, r_{lu}, {\color{red}  \bf \large +}  \, \delta)
                           + \bar{A_1} (\bar{d}_x / 2 {\color{blue}  \bf \large \, -}  \, s, r_{ru}, {\color{blue}  \bf \large -}  \; \delta) \\
            + & \bar{A_1} (\bar{d}_x / 2 {\color{red}  \bf \large \, +}  \,  s, r_{ll}, \, {\color{blue}  \bf \large -}  \; \delta) \,
                           + \bar{A_1} (\bar{d}_x / 2 {\color{blue}  \bf \large \, -}  \, s, r_{rl}, \, {\color{red}  \bf \large +}  \, \delta) \; ) \, / \, 2 \\
            = \quad & \; \dots \; + D \, s^2 + 0 \, s \, \delta + (L + B) \, \delta^2 \; \dots
   \label{Min_A2hat_s_delta}
   \end{split}
\end{equation}
Because $D$ and $L + B$ are positive and the coefficient of the mixed term ist $0$, this $\hat{A_2} ()$
represents a positive definite quadratic form in the minimizing $s, \delta$ with main axis $s = 0, \delta = 0$. \\  

\begin{remark}
Because $\bar{d}_x$ is used instead of $d_x$ in the $\hat{A_2} ()$ above, has not longer sum degree $2$ (the sum of the degrees of $d_x, r_{lu}, r_{ll}, r_{ru}, r_{rl}, \delta$ in a monom). Before doing the following minimazation terms with a sum degree $> 2$ have to be ommited. This leads to $2$ linear equations for the $s, \delta$ to be varied. Alternatively $\bar{d}_x$ has to be used only in the term $B \, \bar{d}_x$, in the terms $D / 2 \, d_x^2, E \, d \, r_{..}, J \, d \, \delta$ we use $d_x$ of degree $1$.
\end{remark}

The $s_{min}, \delta_{min}$ are determined by the now $2$ equations $\frac{d\hat{A_2} (\dots, \; s, \delta)}{ds} = 0$ and $\frac{d\hat{A_2} (\dots, \; s, \delta)}{d\delta} = 0$. Inserting $s_{min}, \delta_{min}$ in \ref{Min_A2hat_s_delta} above we get the following result, for $A_2$ see \ref{A2_d_4r},
$r_l = r_{lu} + r_{ll} - r_{ru} - r_{rl}, r_c = r_{lu} + r_{rl} - r_{ll} - r_{ru}$:
\begin{equation}
   \begin{split}
   \hat{A_2} & (d, r_{lu}, r_{ll}, r_{ru}, r_{rl}) = A_2 (d, r_{lu}, r_{ll}, r_{ru}, r_{rl}) 
                 - \frac{(K  \, r_c - 2 \, B \, d_y)^2}{16 \, (L + B)}  + \; \dots 
   \label{A2hat_d_4r}
   \end{split}
\end{equation}

The minimal area is obtained for the following $s, \delta$:
\begin{equation}
   s_{min}      = - \frac{E \, r_l}{4 \, D}  \qquad
   \delta_{min} = - \frac{K  \, r_c - 2 \, B \, d_y}{4 \, (L + B)}
   \label{s_delta_min}
\end{equation}
The inclination is positive as shown in the upper half of figure \ref{gseg2par}.

\begin{remark}  \label{rem_min2}
The cubic terms in $\bar{A_1}  (d, r, \delta)$, see \ref{PowSerA1_hat} are converted into quadratic terms in
$\frac{d\hat{A_2} (\dots, \; s, \delta)}{ds}$ and $\frac{d\hat{A_2} (\dots, \; s, \delta)}{d\delta}$.
They also appear as quadratic terms in $s_{min}, \delta_{min}$. E.g. an additional term  $d  \, r \,  \delta$ in $\bar{A_1}  (d, r, \delta)$ results as quadratic term $K \, r_c \, (r_{lu} + r_{ru} - r_{ll} - r_{rl}) / (16 \, D \, (L + B))$ in $s_{min}$. But $\hat{A_2} (\bar{d}_x, r_{lu}, r_{ll}, r_{ru},  r_{rl}, \enspace s, \delta)$, see \ref{Min_A2hat_s_delta} has a special form with $\pm s$ and $\pm \delta$ on the rhs.
Therefore linear terms in $s$ and $\delta$ do not occur and so quadratic terms in  $s_{min}, \delta_{min}$ do not result in quadratic terms in the minimum $\hat{A_2} (d, r_{lu}, r_{ll}, r_{ru}, r_{rl})$.
So we can forget about the $10$ cubic terms with 3-fold derivations $\frac{d^3\bar{A_1}}{dd^3} \, d^3 / 6$,
$\frac{d^3\bar{A_1}}{dd^2 \, dr} \, d^2 \, r / 2$, $\dots, \frac{d^3\bar{A_1}}{d\delta^3} \, \delta^3 / 6$, which simplifies a lot!
\end{remark}

\section{2-parameter minimizing, the area of 2 opposite segments \label{2par_area_opp2}}

For the following see also remark \ref{rem_min2}: \\

To include the influence of a second parameter for minimizing, an inclination of the cutted stripe, from the minimal area \ref{A2_d_4r} an additional positive term has to be subtracted, this always lowers the cutted area:
\begin{equation}
   \begin{split}
   \hat{A_2} = & A_2 (d, r_{lu}, r_{ll}, r_{ru}, r_{rl}) - (0.0170374276 \, r_c + 0.2573167207 \, d_y)^2
   \label{AddMin2num}
   \end{split}
\end{equation}

Remark: The area $A_2 ()$ is invariant under permutations of the radii $r_{lu} \leftrightarrow r_{ll}, r_{ru} \leftrightarrow r_{rl}$ on the left and/or right side. The additional term introduced here, is \emph{only} invariant applying these permutations on \emph{both sides simultaneously}. \\
A plausibility check: for upper/lower symmetry, i.e. $r_{ll} = r_{lu}$ and $r_{rl} = r_{ru}$ the additional
term has to be $0$, because the minimum is already received with 1-parameter minimizing. \\

There are 3 types of opposite segments/caps represened by the edges $E_0, E_1, E_2$ in the hexagonal/triangle graph, see \ref{lattice_orientations}. \\

\bibliographystyle{amsplain}

\begin{thebibliography}{10}
\bibitem{CmCsMaVaZs} G. Ambrus, A. Csisz\'{a}rik, M. Matolcsi, D. Varga and P. Zs\'{a}mboki \textit{The Density of planar Sets avoiding Unit Distances}. arXiv:2207.14179 [math.MG]  20 Oct 2022,
 \href {https://doi.org/10.48550/arXiv.2207.14179} {\path{DOI:10.48550/arXiv.2207.14179}}
\bibitem{Cr67} H.T. Croft \textit{Incidence incidents}. Eureka 30, 22-26 (1967)
\bibitem{KaSw2019} B. Kawohl and G. Sweers \textit{On a formula for sets of constant width in 2D}.
 Commun. Pure Appl. Anal. 18,  2117-2131 (2019),
 \href {https://doi.org/10.3934/cpaa.2019095} {\path{DOI:10.3934/cpaa.2019095}},
 \href {http://www.mi.uni-koeln.de/mi/Forschung/Kawohl/kawohl/Gleichdick16.pdf}
       {\path{PDF_at_the_author's_website}}
\end{thebibliography}

\end{document}